\documentclass[a4paper,12pt]{article}

\usepackage{amsfonts}
\usepackage{amsthm}
\usepackage{amsmath}
\usepackage{amsfonts}
\usepackage{latexsym}
\usepackage{amssymb}

\newtheorem{teo}{Theorem}[subsection]
\newtheorem*{teo*}{Theorem}
\newtheorem{lem}[teo]{Lemma}
\newtheorem{cor}[teo]{Corollary}
\newtheorem{pro}[teo]{Proposition}
\theoremstyle{definition}
\newtheorem{fed}[teo]{Definition}
\newtheorem{rem}[teo]{Remark}

\def\noi{\noindent}
\def\QED{\hfill $\blacksquare$}
\def\EOE{\hfill $\blacktriangle$}

\def\eps{\varepsilon}

\def\la{\lambda}
\def\cF{\mathcal{F}}
\def\ga{\gamma}

\def\si{\sigma}

\def\bdem{\begin{proof}}
\def\edem{\renewcommand{\qed}{\hfill $\blacksquare$}
\end{proof}}

\def\N{\mathbb{N}}

\def\R{\mathbb{R}}
\def\C{\mathbb{C}}
\def\inc{\subseteq}
\def\bm{\left(\begin{array}}
\def\em{\end{array}\right)}

\def\eme{\mathcal{M}}

\def\cB{\mathcal{B}}

\def\cD{\mathcal{D}}
\def\cE{\mathcal{E}}

\def\cM{\mathcal{M}}
\def\cN{\mathcal{N}}

\def\cU{\mathcal{U}}
\def\cV{\mathcal{V}}

\def\ese{\mathcal{S}}

\def\ewe{\mathcal{W}}
\def\eme{\mathcal{M}}

\def\ben{\begin{enumerate}}
\def\een{\end{enumerate}}
\def\beq{\begin{equation}}
\def\eeq{\end{equation}}
\def\barr{\begin{array}}
\def\earr{\end{array}}
\def\inv{^{-1}}

\newcommand{\peso}[1]{ \quad \text{ #1 } \quad }
\newcommand{\sub}[2]{{#1}_{\mbox{\tiny{${#2}$}}}}

\DeclareMathOperator{\Preal}{\R\mbox{e}} 

\DeclareMathOperator*{\dist}{dist}

\DeclareMathOperator{\tr}{tr}

\DeclareMathOperator{\sgn}{sgn}

\DeclareMathOperator{\leqp}{\leqslant}
    \DeclareMathOperator{\geqp}{\geqslant}

\def\ga{\gamma}

\newcommand{\pint}[1]{\displaystyle \left \langle #1 \right\rangle}

\newcommand{\hil}{\mathcal{H}}

\newcommand{\cene}{\mathbb{C}^r}
\newcommand{\mat}{\mathcal{M}_r (\C) }
\newcommand{\matu}{\mathcal{U}(r)}
\newcommand{\matsa}{\mathcal{M}_r^{h}(\C)  }
\newcommand{\matah}{\mathcal{M}_r^{ah}(\C)  }

\newcommand{\matinv}{\mathcal{G}\textit{l}\,_r(\C) }

\newcommand{\spec}[1]{\sigma\left( #1\right)}

\newcommand{\conv}{\xrightarrow[n\rightarrow\infty]{}}

\newcommand{\alu}[1]{\Delta\left(#1\right)}
\newcommand{\aluf}[1]{\left|#1\right|^{1/2}U\left|#1\right|^{1/2}}
\newcommand{\aluit}[2]{\Delta^{#1}\left(#2\right)}

\newcommand{\orb}[1]{\ese \left({#1}\right)}
\newcommand{\orbu}[1]{\cU \left({#1}\right)}
\newcommand{\der}[3]{\sub{T}{#2}{#1} \left({#3}\right)}
\newcommand{\dersin}[2]{\sub{T}{#2}{#1}}
\newcommand{\had}[1]{\Psi_{#1}}
\newcommand{\preal}{\sub{P}{\R e}}
\newcommand{\pim}{\sub{P}{\mathbb{I}\mbox{m}}}

\newcommand{\qd}{\sub{Q}{D}}
\newcommand{\qn}{\sub{Q}{N}}
\newcommand{\diagdis}{\cD_{r}^*(\C)}
\newcommand{\bola}[1]{\cB(#1)}

\begin{document}

\title{
\textbf{Convergence of iterated Aluthge transform sequence for diagonalizable matrices} 
\author{Jorge Antezana\thanks{Partially supported by CONICET (PIP 4463/96), Universidad de La PLata (UNLP 11 X350) and
ANPCYT (PICT03-09521).}\and Enrique R. Pujals \thanks {Partially supported by CNPq} \and Demetrio Stojanoff \thanks {Partially supported by CONICET (PIP 4463/96), Universidad de La PLata (UNLP 11 X350) and
ANPCYT (PICT03-09521).}}
}
\maketitle

\noi {\bf Jorge Antezana and Demetrio Stojanoff}

\noi Depto. de Matem\'atica, FCE-UNLP,  La Plata, Argentina
and IAM-CONICET  

\noi e-mail: antezana@mate.unlp.edu.ar and demetrio@mate.unlp.edu.ar

\medskip

\noi {\bf Enrique R. Pujals}

\noi Instituto Nacional de Matem\'atica Pura y Aplicada (IMPA), Rio de Janeiro, Brasil.

\noi e-mail: enrique@impa.br

\vglue1.5truecm

\begin{abstract}
Given an $r\times r$ complex matrix $T$, if $T=U|T|$ is the polar
decomposition of $T$, then, the Aluthge transform is defined by 
$$
\Delta\left(T \right)= |T|^{1/2} U |T |^{1/2}.
$$ 
Let $\Delta^{n}(T)$ denote the n-times iterated Aluthge transform of $T$, i.e.
$\Delta^{0}(T)=T$ and $\Delta^{n}(T)=\Delta(\Delta^{n-1}(T))$, $n\in\mathbb{N}$. 
We prove that the sequence $\{\Delta^{n}(T)\}_{n\in\mathbb{N}}$ converges 
for every $r\times r$ {\bf diagonalizable}  matrix $T$.
We show that the limit $\Delta^{\infty}( \cdot)$ is a map of class $C^\infty$ on the similarity orbit of a diagonalizable matrix, and 
on the (open and dense) 
set of $r\times r$ matrices  with $r$ different eigenvalues. 
\end{abstract}

\vglue2truecm

\noi

\medskip
\noi



\section{Introduction}

Let $\hil$ be a Hilbert space and $T$ a bounded operator defined on $\hil$ whose polar decomposition is $T=U|T|$. The \textit{Aluthge transform} of $T$ is the operator $\alu{T}=|T|^{1/2}U\ |T|^{1/2}$. This was first studied in \cite{[Aluthge]} in relation with the so-called p-hyponormal and log-hyponormal operators. Roughly speaking, the Aluthge transform of an operator is closer to being  normal.

The Aluthge transform has received much attention in recent years. One reason is the connection of Aluthge transform with the invariant subspace problem. Jung, Ko and Pearcy proved in \cite{[JKP0]} that $T$ has a nontrivial invariant subspace if an only if $\alu{T}$ does. On the other hand, Dykema and Schultz proved in \cite{[Dykema]} that the Brown measures is unchanged by the Aluthge transform.

Another reason is related with the iterated Aluthge transform. 
Let $\aluit{0}{T}=T$ and $\aluit{n}{T}=\alu{\aluit{n-1}{T}}$ 
for every $n\in\N$. 
%
It was conjectured in \cite{[JKP0]} that the sequence 
$\{\aluit{n}{T}\}_{n\in\N}$ converge in the norm topology. 
Although this conjecture was stated for operators on an arbitrary Hilbert space, 
it was corrected and restated for matrices in \cite{[JKP1]} by Jung Ko and 
Pearcy and receantly extended to finite factors in \cite{[Dykema]} by Dykema and Schultz. 
In these spaces, it still remains open and there only exist some partial results. 
For instance, Ando and Yamazaki proved in \cite{[Ando-Yamaza]} that the conjecture 
is true for $2\times 2$ matrices and Dykema and Schultz in \cite{[Dykema]} proved 
that the conjecture is true for an operator $T$ in a finite factor such that 
the unitary part of its polar decomposition normalizes an abelian subalgebra 
that contains $|T|$. (see \cite{[Ando]}, \cite{[Wu]}  and \cite{[Yamazaki]} 
for other results that support the conjecture in finite factors).


A result proved independently by Jung, Ko and Pearcy in \cite{[JKP1]}, and  by Ando in \cite{[Ando]},
states that, given an $r\times r$ matrix $T$, the limit 
points of the sequence $\{\aluit{n}{T}\}_{n\in\N}$ are normal matrices with the same 
characteristic polynomial as $T$. In particular, if the sequence of iterated Aluthge 
transforms converge, the limit function, defined by $\displaystyle T\mapsto 
\lim_{n\to\infty}\aluit{n}{T}$, whould be a retraction from the space of matrices 
onto the set of normal operators.

\medskip
Another important result, concerning the finite dimensional case, states that it is enough 
to prove the conjecture for invertible matrices (see for example \cite{[AMS]}). 
Note that, for an invertible matrix $T$
$$
\alu{T}=|T|^{1/2} \, T\, |T|^{-1/2}.
$$ 
So the Aluthge transform of $T$ belongs to the similarity orbit of $T$. This suggest 
that we can study the Aluthge transform restricted to the similarity orbit of some invertible operator. 

From that point of view, the diagonalizable case has some advantages. First of all, note that the 
similarity orbit of a diagonalizable operator contains a compact submanifold of fixed points, and 
the sequence $\{\aluit{n}{T}\}_{n\in\N}$ goes to this submanifold as $n\to\infty$. In fact, since 
$T$ is diagonalizable, the similarity orbit of $T$ coincides with the similarity orbit of some 
diagonal operator $D$, which we denote $\orb{D}$. The unitary orbit of $D$, denoted by $\orbu{D}$, 
is a compact submanifold of $\orb{D}$ that consists of all normal matrices in $\orb{D}$. Hence 
$\orbu{D}$ is  fixed by the Aluthge transform and the limits points of 
the sequence $\{\aluit{n}{T}\}_{n\in\N}$ belongs to $\orbu{D}$. In contrast, for 
non-diagonalizable operators, the similarity orbit 
does not have fixed points, and the sequence of iterated Aluthge transforms goes to points 
that do not belong to the similarity orbit.

On the other hand, numerical computations, as well as Ando-Yamazaki's $2\times 2$ computations
(see \cite{[Ando-Yamaza]}), suggest that the rate of convergence of the sequence $\{\aluit{n}{T}\}_{n\in\N}$, 
for diagonalizable operators $T$, becomes exponential after some iterations. However, 
it seems that this behavior is not shared by the non-diagonalizable case.

For these reasons, we decided to study the diagonalizable case. Note that if we restrict the Aluthge transform to the similarity orbit of an invertible diagonalizable matrix $T$, a dynamical system approach can be performed.

In fact, we show that for any $N\in \orbu{D}$ there is a local submanifold 
$\ewe_N^s$ transversal to $\orbu{D}$ characterized by the matrices that 
converges with a exponential rate to $N$ by the iteration of the 
Aluthge transform. Moreover, the union of these submanifolds form an open 
neighbourhood of $\orbu{D}$ (see Corollary \ref{entorno y algo mas}). 
Thus, since the sequence $\{\aluit{n}{T}\}_{n\in\N}$ goes towards $\orbu{D}$, 
for some $n_0$ large enough the sequence of iterated Aluthge tranforms enters 
this open neighborhood and converge exponentially. 

These results follow from the classical arguments of stable manifolds (first introduced independently  
by Hadamard and Perron, see theorem \ref{pseudohyperbolic}; for details and 
general results 
about the stable manifold theorem see \cite{[HPS]} or the Appendix at the end 
of this work). 
To conclude that, it is shown that the 
derivative of the  Aluthge transform in any $N\in\orbu{D}$ has two invariant complementary 
directions, one tangent to $\orbu{D}$, and other transversal to it, where the derivative 
is a contraction (see Theorem \ref{the key}). Using these results, 
we prove that the sequence $\{\Delta^{n}(T)\}_{n\in\mathbb{N}}$ converges 
for every $r\times r$ {\bf diagonalizable}  matrix $T$.
We also show that the limit $\Delta^{\infty}( \cdot)$ is a map of class $C^\infty$ 
on the similarity orbit of a diagonalizable matrix, and 
on the (open and dense) set of $r\times r$ matrices  with $r$ different 
eigenvalues. 

This paper is organized as follows: 
in section 2, we collect several preliminary definitions and results about the
the stable manifold theorem, about the geometry of similarity and unitary orbits, 
and about known results on Aluthge transform. 
In section 3, we prove the convergence results and we study the smoothness of 
the limit map $T \mapsto \aluit{\infty}{T}$, mainly for $r\times r$ matrices 
with $r$ different eigenvalues. 
The basic tool, to apply the stable manifold theorem to the 
similarity orbit of a diagonal matrix, is the mentioned  
Theorem \ref{the key}, whose proof, somewhat technical, is done in 
section 4. In the Appendix, 
we sketch the proof of the classical version of the stable 
manifold theorem in order to show how it can be modified 
in our context, where the invariant set is a smooth submanifold 
consisting of fixed points, getting stronger results 
on the regularity conditions of the prelamination $\{\ewe_N^s\}_{N\in \orbu{D}}$.
  
We would like to aknowledge Prof. M. Shub for comments 
and suggestion about the stable manifold theorems, and 
Prof. G. Corach who told us  about the Aluthge transform, and 
shared with us fruitful discussions concerning these matters.

\section{Preliminaries.}
In this paper $\mat$ denotes the algebra of complex $r\times r$ matrices, 
$\matinv$ the group of all invertible elements of $\mat$, $\matu$ the group 
of unitary operators, and $\matsa$ (resp. $\matah$) denotes the real algebra of hermitian (resp. antihermitian) matrices. 
Given $T \in \mat$, $R(T)$ denotes the
range or image of $T$, $\ker(T)$ the null space of $T$, 
$\sigma (T)$ the spectrum of $T$, $\tr(T)$ the trace of $T$,
and $T^*$ the adjoint of $T$. If $v \in \cene$, we debote by 
$\mbox{\rm diag}(v) \in \mat$ the diagonal matrix with $v$ in its diagonal. 
We shall consider the space of matrices $\mat$ as a real Hilbert space 
with the inner product defined by
$$
\pint{A,\ B}=\Preal\big(\tr(B^*A)\big).
$$
The norm induced by this inner product is the so-called Frobenius norm, that is denoted by $\|\cdot \|_2$.
Along this note we also use the fact that every subspace $ \ese $ of $ \mathbb C^n$
induces a representation of elements of $\mat$ by $2 \times 2$ block
matrices, that is, we shall identify each  $A\in\mat$ with
a $2\times 2$-block matrix
$$
\begin{pmatrix}
  A_{11} & A_{12} \\
  A_{21} & A_{22}
\end{pmatrix}\begin{array}{cc}
  \ese  \\
  \ese^\bot
\end{array} ,
$$ 
where $A_{11}=\left. A\right|_{\ese,\,\ese}\,$, $A_{12}=\left. A\right|_{\ese^\bot,\,\ese}\,$, $A_{21}=\left. A\right|_{\ese,\,\ese^\bot}\,$ and $A_{22}=\left. A\right|_{\ese^\bot,\,\ese^\bot}\,$.


\medskip
On the other hand, let $M$ be a manifold. By means of $TM$ we denote the tangent bundle of $M$ and by means of $T_xM$ we denote the tangent space at the point $x\in M$. Given a function $f\in C^{r}(M)$, where $r=1,\ldots,\infty$, $\der{f}{x}{v}$ denotes the derivative of $f$ at the point $x$ applied to the vector $v$.

\subsection{Stable manifold theorem}

In this section we state the stable manifold theorem for an invariant  set of a smooth 
endomorphism (see \ref{teorema 5.5} below). The stable set is naturally defined for a 
fixed point of an endomorphism, as the set of points with positive trajectories heading 
directly towards the fixed point. This notion is the natural extension of the stable 
eigenspaces of a linear transformation (the ones associated to the eigenvectors with 
modulus smaller than one)  into the nonlinear regimen. In fact, 
a natural intuitive approach to the idea of the stable manifold is to consider a fixed point of a smooth differentiable map  such that the derivative of the map at the fixed point has absolute value smaller than one. In this case, the linear map induced by the derivative is a map that share the same fixed point and such that any trajectory converges by forward iterate to the fixed point with a exponential rate of contraction. Using that the linear map is a ``good approximation of the map in a small neighborhood of the fixed point", it  follows that the map has the same dynamical behavior of its linear part.

A more general approach is based in the techniques known as graph transform 
operator. This approach can be naturally extended for invariant sets, being 
almost straightforward when the set consists of fixed points. An sketched 
version of the proof of Theorem \ref{teorema 5.5}, using these techniques, is 
done in the Appendix at the end of this work 
(see also \cite[Thm 5.5]{[HPS]}).



\medskip
\noi
Let $M$ be a smooth Riemann manifold and $N\subseteq M$ a submanifold (not 
necessarily compact). Throughout this subsection $\sub{T}{N}M$ denotes the 
tangent bundle of $M$ restricted to $N$.

\begin{fed}\label{prelamination}
A $C^r$ \textit{pre-lamination} indexed by $N$ is a continuous choice of a $C^r$ embedded disc $\cB_x$ 
through each $x\in N$. Continuity means that $N$ is covered by open sets  $\cU$ 
in which $x\to B_x$ is given by
$$
\cB_x=\sigma(x)((-\eps,\eps)^k)
$$
where $\sigma: \cU \cap N\to \mbox{Emb}^r((-\eps,\eps)^k,M)$ is a continuous section. 
Note that $\mbox{Emb}^r((-\eps,\eps)^k,M)$ is a $C^r$ fiber bundle over $M$ whose 
projection is $\beta\to \beta(0)$. Thus $\sigma(x)(0)=x$. 
If the sections mentioned above are $C^s$, $1\leq s\leq r$, we say that the $C^r$ 
pre-lamination is of class $C^s$.

\end{fed}

\begin{fed}\label{self coherent}
A prelamination is \textit{self coherent} if the interiors of each pair of its discs meet in a relatively open subset of each.
\end{fed}


\begin{fed}\label{pseudohyperbolic}
Let $f$ be a smooth endomorphism of $M$, $\rho>0$, and suppose that $\left. f\right|_{N}$ is a homeomorphism. Then, $N$ is  \textit{$\rho$-pseudo hyperbolic} for $f$ if there exist two smooth subbundles of $\sub{T}{N}M$, denoted by $\cE^s$ and $\cF$, such that 
\begin{enumerate}
	\item $\sub{T}{N}M= \cE^s \oplus \cF$;
        \item $\sub{T}{}N=  \cF$;
	\item Both, $\cE^s$ and $\cF$, are $Tf$-invariant;
	\item $T\, f$ restricted to $\cF$ is an automorphism, which expand it by a factor
	greater than $\rho$. 
	\item $\dersin{f}{x}:\cE_x^s\to\cE_{f(x)}^s$ has norm lower than $\rho$. \EOE
\end{enumerate}
\end{fed}

\noi
Observe that if $N$ is  \textit{$\rho$-pseudo hyperbolic} then there exists a positive constant $\la=\la(\rho)<1$ such that 
\begin{eqnarray}\label{dominacion}
\frac{||Df_{\cE^s}||}{m(Df_{|\cF})}<\la \ ,
\end{eqnarray}
where $m(.)$ means the minimum norm. 
If $N$ consists of fixed points then, for example,  $N$ is  \textit{$\rho$-pseudo hyperbolic} 
(also called \textit{normally hyperbolic}) if  there is a $Tf-$invariant 
subbundle $\cE^s$ (of $T_N M$) 
complement  to $TN$, such  that $Tf$ contracts  more sharply than any contraction in $TN$.  
In the case that $\cE^s$ is uniformly contracted, it follows that for any point $x\in N$ it 
is possible to find an $f-$invariant submanifold transversal to $N$ tangent to $\cE^s$ 
and characterized as the set of points with trajectories asymptotic to the trajectory of $x.$

\begin{teo}[Stable manifold theorem]\label{teorema 5.5}
Let $f$ be a $C^r$ endomorphism of $M$ with a $\rho$-pseudo hyperbolic submanifold 
$N$ with  $\rho< 1$. 
Then, there is a $f$-invariant and self coherent 
$C^r$-pre-lamination of class $C^0$,   $\ewe^s:  N\to \mbox{Emb}^r((-1,1)^k,M)$ 
such that, for every $x \in N$, 
\begin{enumerate}
\item $\ewe^s(x)(0)=x$,
\item  $\ewe_x^s=\ewe^s(x)((-1,1)^k)$  is tangent to $\cE_x^s$  at every $x\in N$,
\item $
\ewe_x^s \inc \Big\{y\in M\ : \  \dist(f^n(x),f^n(y))<\dist(x,y)\rho^n\Big\}
$.
\end{enumerate}

\end{teo}
\proof 
See the proof in subsection \ref{teorema 5.5ap} of the Appendix. \QED

\begin{cor}[Smoothness of the stable lamination for a submanifold of fixed points]\label{teorema 5.5bisbis}
Let $f$, $M$ and $N$ as in Theorem \ref{teorema 5.5}. 
Let us assume that any point $p$ in $N$ is a {\bf fixed point}.  
%
Then $C^r$-pre-lamination $\ewe^s: \cN \to \mbox{Emb}^r((-1,1)^k,M)$ is  of class $C^r$.
\end{cor}

\proof 
See Corollary \ref{teorema 5.5bisbis2} in the Appendix. \QED

\begin{rem}\label{aclaraciones}
Observe that, from Theorem \ref{teorema 5.5}, it holds that, for every $x\in 
	N$ 
	$$
	\sub{T}{x}\ewe_x^s=\cE_x^s \ . 
	$$
	If $N$ consists on fixed pionts, from the regularity conditions of the
	pre-lamination $\{\ewe_x^s\}_{x\in N}$ assured by Corollary \ref{teorema 5.5bisbis}, 
	 we get that, for any  $x\in N$,  there exists $\ga>0$ such that  
	$$B(x,\ga)\subset \bigcup_{x\in N} \ewe_x^s \ .	$$
In other words, it means that $\bigcup_{x\in N} \ewe_x^s$ contains 
an open neighborhood $\ewe(N)$ of $N$ in $M$. 
Therefore,  condition 3 of Theorem \ref{teorema 5.5} implies that 
for every $x\in N$ there exists an open neighborhood $\cU$ 
of $x$ (open relative to $M$) such that
\begin{equation}\label{converge}
	\ewe_x^s\cap \cU=\Big\{y\in \cU:\ \dist(x,f^n(y))< \  \dist(x,y)\,  \rho^n \Big\}.
	\end{equation}
In particular, 
$\ewe_x^s\cap \ewe_y^s = \varnothing  $ if $x \neq y$. Moreover, we can assure that 
the (well defined) map 
\begin{equation}\label{proyecta}
p : \ewe(N)  \to N \peso{ given by }
p(a) = x \peso{if} a \in  \ewe_x^s(x) 
\end{equation}
is of class $C^r$. 
	\EOE
\end{rem}

\subsection{Similarity orbit of a diagonal matrix}

In this subsection we recall some facts about the similarity orbit of a diagonal matrix. 

\begin{fed}
Let  $D\in\mat$. 
By means of $\orb{D}$ we denote the similarity orbit of $D$:  
$$
\orb{D} = \{ \ SDS\inv\  : \ S \in \matinv \ \} \ . 
$$
On the other hand, $\orbu{D}= \{ \ UDU^* \  : \ U \in \matu \ \}$ 
denotes the unitary orbit of $D$. 
We donote by 
$\sub{\pi}{D} : \matinv \to \orb{D} \inc \mat $  the $C^\infty$ map defined by 
$\sub{\pi}{D}(S) = SDS\inv $.  
With the same name we note its resrtiction to the unitary group:  
$\sub{\pi}{D} : \matu \to \orbu{D} $. 
\EOE
\end{fed}

\begin{pro}\label{son variedades}
The similarity orbit $\orb{D}$ is a $C^\infty$ submanifold of $\mat$, and the projection
$\sub{\pi}{D} : \matinv \to \orb{D}$ becomes a submersion. Moreover, $\orbu{D}$ is a compact submanifold of $\orb{D}$, which consists of the normal elements of $\orb{D}$, and 
$\sub{\pi}{D} : \matu \to \orbu{D}$ 
is a submersion.\QED
\end{pro}

\noi
For every $N=UDU^*\in\orbu{D}$, it is well known (and easy to see) that
\begin{align}
\sub{T}{N}\,\orb{D}&= \sub{T}{I}( \pi_N ) (\mat\,)  =\{[A,N]=AN-NA: \ A\in\mat\}. \nonumber\\
\intertext{In particular}
\sub{T}{D}\,\orb{D}&=\{AD-DA: \ A\in\mat\} \nonumber
\\\label{tan1} &=\{X\in\mat:\ X_{ij}=0 \ \mbox{for every $(i,j)$ such that $d_i=d_j$}\}.\\
\intertext{Note that, }
\sub{T}{N}\,\orb{D}&=\{[A,N]=AN-NA: \ A\in\mat\}\nonumber
\\&=\{(UBU^*)UDU^*-UDU^*(UBU^*): \ B\in\mat\}\nonumber
\\&=\{U[B,D]U^*=BD-DB: \ B\in\mat\}
=U\Big(\sub{T}{D}\, \orb{D}\Big)U^*\ .\\
\intertext{On the other hand, since $\sub{T}{I}\, \matu = \matah= \{ A \in \mat : A^* = -A\}\,$, we obtain}
\sub{T}{D}\,\orbu{D}
&= \sub{T}{I}( \pi_D ) (\matah \,) =\{[A,D]=AD-DA: \ A\in\matah\} \peso{and} 
\ ,\nonumber\\
\label{ss}\sub{T}{N}\,\orbu{D}&=\{[A,N]=AN-NA: \ A\in\matah\}=U\Big(\sub{T}{D}\,\orbu{D}\Big) U^* \ . 
\end{align}

\noi Finally, along this paper we shall consider on $\orb{D}$ (and in $\orbu{D}$) the Riemannian structure inherited from $\mat$ (using the usual inner product on their tangent spaces). For $S, T \in \orb{D}$, we denote by $\dist(S, T)$  the Riemannian distance between $S$ and $T$ (in $\orb{D}\,$). 
Observe that, for every $U \in \matu$, one has that $U\orb{D} U^* = \orb{D} $ and the map 
$T \mapsto UTU^*$ is isometric, on $\orb{D}$,  with  respect to the Riemannian metric as well as with respect to the $\| \cdot \|_2$ metric of $\mat$. 

\subsection{Definition and basic facts about Aluthge transforms}

\begin{fed}\rm
Let $T\in\mat$, and suppose that $T=U|T|$ is the polar
decomposition of $T$. Then, we define the
Aluthge transform of $T$ in the following way:
\begin{align*}
\alu{T}&=\aluf{T}
\end{align*}
On the other hand, $\aluit{n}{T}$ denotes the n-times iterated Aluthge transform of $T$, i.e.
\begin{align*}
\aluit{0}{T}=T; \peso{and}
\aluit{n}{T}=\alu{\aluit{n-1}{T}}\quad n\in\N.
\end{align*}
\end{fed}

\noi The following proposition contains some properties of Aluthge transforms which follows easily from its definition.

\begin{pro}\label{facilongas}
Let $T\in\mat$. Then:
\begin{enumerate}
  \item $\alu{cT}=c\alu{T}$ for every $c\in \C$.
  \item $\alu{VTV^*}=V\alu{T}V^*$ for every $V\in\matu$.
  \item If $T=T_1\oplus T_2$ then    
  $\alu{T}=\alu{T_1}\oplus \alu{T_2}$.
  \item $\|\alu{T}\|_2\leqp \|T\|_2$.
  \item $T$ and $\alu{T}$ have the same  characteristic polynomial, in particular,
  $\spec{\alu{T}}=\spec{T}$.
\end{enumerate}
\end{pro}

\noi The following theorem states the regularity properties of Aluthge transforms (see \cite{[Dykema]}).

\begin{teo}\label{continuidad}
The Aluthge transform is $(\|\cdot\|_2\ ,\ \|\cdot\|_2)$-continuous in $\mat$ and it is of class $C^\infty$ in $\matinv$.
\end{teo}

\noi Now, we recall a result proved independently by Jung, Ko and Pearcy in \cite{[JKP1]}, and  by Ando in \cite{[Ando]}.

\begin{pro}\label{puntos limites normales}
If $T\in\mat$,  the limit points of the sequence
$\{\aluit{n}{T}\}_{n\in \N}$ are normal. Moreover, if $L$ is a limit point,
then $\spec{L}=\spec{T}$ with the same algebraic multiplicity.
\end{pro}

\noi Finally, we mention a result concerning the Jordan structure of Aluthge transforms proved in \cite{[AMS]}. We need the following definitions. 

\begin{fed}\rm Let $T \in \mat $ and $\mu \in \spec{T}$. We denote
\ben
\item $m(T, \mu)$ the $algebraic$ $multiplicity$ of the eigenvalue $\mu$ for $T$.
\item $m_0(T, \mu)= \dim \ker( T-\mu I )$,
the $geometric$ $multiplicity$ of $\mu$.
\een
\end{fed}

\begin{pro}\label{algeo}
Let $T\in \mat$. 
\ben
\item If $\ 0 \in \spec{T}$, then, there exists $n\in\N$ such that
$$
 m(T, 0)=  m_0(\aluit{n}{T}, 0) = \dim \ker(  \aluit{n}{T}).
$$ 
\item For every  $\mu \in \sigma(T)$, 
$m_0(T, \mu ) \leqp m_0( \alu{T}, \mu ).$  
\een
\end{pro}

\noi
Observe that this implies that, if $T$ is diagonalizable (i.e. $m_0(T, \mu )  = m(T, \mu ) $ for every $\mu$), 
then also $\alu{T}$ is diagonalizable.

\section{The iterated Aluthge transform} 

\subsection{Convergence of iterated Aluthge transform sequence for diagonalizable matrices}\label{convergencia}

\noi In this section, we prove the convergence of iterated Aluthge transforms for diagonalizable matrices. The key tool, which allows to use the stable manifold theorem
\ref{teorema 5.5}, is the following theorem, whose proof is rather long and technical. 
For this reason, we postpone it until section \ref{la prueba}, and we continue 
in this section with its consequences.

\begin{teo}\label{the key}
Let $D= \mbox{\rm diag}(d_1,\ldots,d_r)  \in \mat$ be an invertible diagonal matrix.
The Aluthge transform $\alu{\cdot}:\orb{D}\to\orb{D}$ 
is a $C^\infty$ map. For every $N\in\orbu{D}$, there exists a 
subspace $\sub{\cE}{N}^s$ of the tangent space $\sub{T}{N}\orb{D}$ such that
\begin{enumerate}
	\item $\sub{T}{N}\orb{D}=\sub{\cE}{N}^s\oplus \sub{T}{N}\orbu{D}$;
	\item Both, $\sub{\cE}{N}^s$ and  $\sub{T}{N}\orbu{D}$, are $T\,\Delta$-invariant;
	\item $\left\|\left. T\,\Delta\right|_{\sub{\cE}{N}^s}\right\|\leq\sub{k}{D}<1$, where $\displaystyle \sub{k}{D}=\max_{i,\,j\,:\ d_i\neq d_j}\frac{|1+e^{i(\arg(d_j)-\arg(d_i))}|\, |d_i|^{1/2}|d_j|^{1/2}}{|d_i|+|d_j|}$;
  \item If $U\in\matu$ satisfies $N=UDU^*$, then $\sub{\cE}{N}^s=U(\sub{\cE}{D}^s)U^*$.
  \end{enumerate}
  In particular, the map $\orbu{D} \ni N \mapsto \sub{\cE}{N}^s$ is smooth. This fact can be 
  formulated in terms of the projections $P_N$ onto $\sub{\cE}{N}^s$ parallel to $\sub{T}{N}\orbu{D}$, 
  $N \in \orbu{D}$. \QED
\end{teo}
\begin{cor}\label{entorno y algo mas}
Let $D= \mbox{\rm diag}(d_1,\ldots,d_r)  \in \mat$ be an invertible diagonal matrix.
Let $\sub{\cE}{N}^s$ and $k_D$ as in Theorem \ref{the key}. Then, in 
$\orb{D}$ there exists a $\Delta$-invariant $C^\infty$-pre-lamination $\{\ewe_N\}_{N\in\orbu{D}}$ of class $C^\infty$ such that, for every $N \in \orbu{D}$, 
\begin{enumerate}
	\item $\ewe_N$ is a $C^\infty$ submanifold of $\orb{D}$.
	\item $\sub{T}{N}\ewe_N=\sub{\cE}{N}^s\,$.
	\item   If \  $k_D <\rho<1$, then $\dist(\aluit{n}{T}-N)\leq \dist(T, N)  \rho^n$, for 
	every $T\in \ewe_N\,$.
	 
	\item If $N_1\neq N_2$ then $\ewe_{N_1}\cap\ewe_{N_2}=\varnothing$.
  \item There exists an open subset $\ewe(D)$ of $\orb{D}$ such that 
  \ben
  \item [a. ]
  $\orbu{D} \inc  \ewe(D) \inc \displaystyle \bigcup_{N\in\orbu{D}}\ewe_N$, and 
  \item [b. ] The projection $p:\ewe(D)\to \orbu{D}$, defined by 
  $p(T)=N$ if $\,T\in\ewe_N$, is of 
  class $C^\infty$.
  \een
\end{enumerate}
\end{cor}
\proof By Theorem \ref{the key}, for every $k_D< \rho < 1$, 
$\orbu{D}$ is $\rho$-pseudo hyperbolic for $\Delta$
(see Definition \ref{pseudohyperbolic}), and it consists of fixed points. Thus, by 
Corollary  \ref{teorema 5.5bisbis} and Remark \ref{aclaraciones}, 
we get a $C^\infty$ and $\Delta$-invariant prelamination of class $C^\infty$, 
$\{\ewe_N\}_{N\in\orbu{D}}$ which satisfies all the properties 
of our statement.  
\QED

\bigskip
\noi
In order to prove the convergence of iterated Aluthge transforms for diagonalizable matrices, we first reduce the problem to the invertible case.  
In \cite{[AMS]} it was proved that if the sequence of iterated Aluthge transforms converge 
for every invertible matrix, then it 
converge for every matrix. In our case, we need to prove that if the sequence of iterated 
Aluthge transforms converge for every diagonalizable invertible matrix, then it does for 
every diagonalizable matrix. The proof of the second statement is essentially the same 
as the previous one, but, for a sake of completeness, we include its proof.

\begin{lem}\label{reeduccion de conjetura}
 If the sequence $\{\aluit{m}{S}\}_{m\in\N}$
 converges for every diagonalizable invertible matrix
 $S\in  \mat $ and every $r \in \N$, then the
 sequence $\{\aluit{m}{T}\}_{m\in\N}$ converges for every diagonalizable matrices $T\in \mat $ 
 and every $r \in \N$.
\end{lem}
\bdem
Let $T \in \mat$. As we have observed after 
Proposition \ref{algeo}, 
if $T$ is diagonalizable, then $\alu{T}$ is also diagonalizable. So, if we begin 
with a diagonalizable matrix $T$, then every element of the sequence 
$\{\aluit{m}{T}\}_{m\in\N}$ is diagonalizable. By Proposition \ref{algeo}, 
we can also assume that $m(T, 0) = m_0(T, 0)$.
Note that, in this case,
$\ker( \alu{T})=\ker (T)$ because $\ker (T)\subseteq \ker (\alu{T})$
and $m(\alu{T}, 0 ) = m(T, 0 )$.
On the other hand, $R(\alu{T})\subseteq R(|T|)$ so that
$R(\alu{T})$ and $\ker(\alu{T})$ are orthogonal subspaces. Thus,
there exists a unitary matrix $U$ such that 
$$ 
U\alu{T} U^*=\begin{pmatrix}S &0 \\ 0&0\end{pmatrix} 
$$ 
where $S\in M_s(\C )$ is invertible and diagonalizable ($s = n - m(T, 0)$\;). Since  for every $m \geqp 2$ 
$$ 
\aluit{m}{T}=U^* \begin{pmatrix}\aluit{m-1}{S} &0 \\
0&0\end{pmatrix}U \ , 
$$ 
the sequence $\{\aluit{m}{T}\}$ converges,
because the sequence $\{\aluit{m-1}{S}\}$ converges by 
hypothesis.\edem

\medskip
\begin{teo}\label{convergencia para diagonalizables}
Let $T\in\mat$ be a {\bf diagonalizable} matrix. Then $\{\aluit{n}{T}\}_{n\in\N}$ converges.
\end{teo}
\bdem
Using Lemma \ref{reeduccion de conjetura},  we can assume that $T$ is invertible. 
Then, $T\in \orb{D}$ for some invertible diagonal matrix $D$.
By Corollary  \ref{entorno y algo mas} and Remark \ref{aclaraciones}, 
we get on $\orb{D}$ a $C^\infty$ and $\Delta$-invariant prelamination 
of class $C^\infty$, denoted by $\{\ewe_N\}_{N\in\orbu{D}}$, such that 
\ben
\item The set $\bigcup _{N \in \orbu{D}} \ewe_N $ contains an open neighborhood
$\ewe(D)$ of $\orbu{D}$ in $\orb{D}$.
\item If \  $k_D <\rho<1$, then $\|\aluit{n}{A} - N \|_2 \le
\dist(\aluit{n}{A}-N)\leq \dist(A, N)  \rho^n$, for 
	every $A\in \ewe_N\,$.
\een
On the other hand, by Proposition \ref{puntos limites normales}, 
there exists $m \in \N$ such that  
$A = \aluit{m}{T} \in \bigcup _{N \in \orbu{D}} \ewe_N $.
Thus, for $n >m$, $\aluit{n}{T} = \aluit{n-m}{A} \conv N$, where 
$N \in \orbu{D}$ is  the unique  element of $\orbu {D}$ such that $A \in \ewe_N$. 
\edem

\begin{rem}
From Theorem \ref{convergencia para diagonalizables} 
it can be deduced Ando and Yamazaki's result on the convergence of the 
iterated Aluthge sequence for $2\times 2$ matrices. Indeed, in $\eme_2(\C)$, the spectrum of 
matrices uncovered by Theorem \ref{convergencia para diagonalizables} must be a singleton. 
Therefore, by Proposition \ref{puntos limites normales}, the iterated Aluthge sequence for 
those matrices has only one limit point. So, it converges.\EOE
\end{rem}





\begin{pro}\label{retraccion}
Let $D\in \mat$ be diagonal and invertible. 
Then the sequence $\{\Delta ^n\}_{n\in\N}$, 
resticted to the similarity orbit $\orb{D}$,  converges uniformly on 
compact sets to a $C^\infty$ limit function $\Delta^\infty:\orb{D}\to\orbu{D}$. 
In particular, $\Delta^\infty$ is a $C^\infty$ retraction from $\orb{D}$ onto $\orbu{D}$.
\end{pro}
\bdem
Let $\Delta^\infty$ 
be the limit function, which exists by Theorem \ref{convergencia para diagonalizables}.
We can apply Corollary \ref{entorno y algo mas}, and we shall use 
its notations. Fix $T\in\orb{D}$. 
By Proposition \ref{puntos limites normales} there 
exists $k\in\N$ such that $\aluit{k}{T}\in\ewe(D)$. 
  By the continuity of $\alu{\cdot}$, there exists a 
neighborhood $\cU$ of $T$ such that $\aluit{k}{\cU}\subseteq \ewe(D)$. Hence, 
if $p$ is the projection defined in Corollary \ref{entorno y algo mas}, $\left.\Delta^\infty\right|_{\cU}=\left.(p\circ \Delta^k)\right|_\cU\ $, 
which proves that the map $\Delta^\infty$ is  $C^\infty$ at $T$.

On the other hand, to prove that the convergence of $\{\aluit{n}{\cdot}\}_{n\in\N}$ 
is uniform on compact sets, suppose that $\cU$ has compact closure, and denote by 
$$C = \sup \{ \dist( \aluit{k }{S}, \aluit{\infty }{S}): S \in \cU \} \ .
$$ 
Fix $\eps>0$ and take $m_0>k$ such that 
$Ck_D^{m_0-k}<\eps$. Then, using (4) of Corollary \ref{entorno y algo mas},
for every $m\geq m_0$ and every $S\in\cU$
$$
\dist(\aluit{m}{S}-\aluit{\infty}{S})= \dist\big(\aluit{m-k }{\aluit{k }{S}\,}-\aluit{\infty}{\aluit{k }{S}\,}\big) \leq \eps.
$$
This proves that for every $T\in\orb{D}$ there exists a neighborhood of $T$ where the 
convergence is uniform. Therefore, by standard arguments, it follows that the convergence 
is uniform on compact sets. \edem

\begin{rem}\rm Let $D\in \mat$ be diagonal but not invertible. If $T\in \orb{D}$, 
by arguments similar to those used in the proofs of Lemma \ref{reeduccion de conjetura}
and Proposition \ref{retraccion} 
it can be proved that $\alu{T}\in \orb{D}$, and the map 
$\Delta^{\infty} \big|_{\orb{D}} : \orb{D} \to \orbu{D}$ 
is a retraction of calss $C^\infty$. 
\end{rem}

\subsection{Smoothness  of the map $T \mapsto \aluit{\infty}{T}$ on $\diagdis$}

Let  $\diagdis$ be the set of diagonalizable and invertible matrices in $\mat$ 
with $r$ different eigenvalues (i.e. every eigenvalue has algebraic multiplicity equal to one). Observe that $\diagdis$ is an open dense subset of $\mat$ and it is invariant by the  Aluthge transform. 
If $\aluit{\infty}{\cdot}$  denotes the limit of the sequence of iterated 
Aluthge transforms, which is defined on the set of diagonalizable matrices by 
Theorem \ref{convergencia para diagonalizables}, we shall show that 
$T \mapsto \aluit{\infty}{T}$ is of class $C^\infty$ on $\diagdis$.
The proof of this result  essentially follows the same lines as 
Proposition \ref{retraccion}. For this reason, we expose a sketched version of the proof, where we
only point out the main differences.

We already know that the map $\aluit{\infty}{\cdot}$ is of class $C^\infty$
if it is restricted to the orbits $\orb{T}$ for any  $T \in \diagdis$. 
In order to study the behavior of this map outside the orbit of $T$, we need to define 
the following sets: let $D\in \diagdis$ be a diagonal matrix and let $\eps >0$;  then 
\begin{align*}
\bola{D,\, \eps}&=\Big\{D'\in \diagdis : \ D'\ 
\mbox{is diagonal and $\|D-D'\|_2<\eps$}\Big\};\\
\orb{D,\, \eps}&= \Big\{SD'S^{-1}:\ D'\in \bola{D,\,\eps}\ \mbox{and}\ S\in\matinv\Big\}
= \bigcup _{ D'\in \bola{D,\,\eps}} \orb{D'}  ;\\
\orbu{D,\, \eps}&=\Big\{UD'U^*:\ D'\in\bola{D,\,\eps}\ 
\mbox{and}\ U\in\matu\Big\} = \bigcup _{ D'\in \bola{D,\,\eps}} \orbu{D'} \ .
\end{align*}
The set  $\orb{D,\, \eps}$ is  invariant for $\alu{\cdot}$ and it is also open in $\matinv$ for $\eps$ small enough. 
Since $D \in \diagdis$, it can be proved that $\orbu{D,\, \eps}$ is a 
smooth submanifold of $\mat$, and it consists on the fixed points of  $\orb{D,\, \eps}$. 
For each $N \in \orbu{D,\, \eps}$, if $\{N\}'$ denotes the 
subspace $\{A \in \mat : AN = NA\}$,  the tangent space  
$\sub{T}{N}\orbu{D,\,\eps}$ can be decomposed as
$
\sub{T}{N}\orbu{D,\,\eps} = \sub{T}{N}\orbu{D} \oplus  \{N\}' \ .
$
Then, $\sub{T}{N}\orb{D,\,\eps}=\mat$ can be decomposed as 
\begin{equation}\label{summma}
\sub{T}{N}\orb{D,\,\eps} = \sub{T}{N}\orb{D} \oplus \{N\}'
= \Big( \sub{\cE}{N}^s\oplus \sub{T}{N}\orbu{D}\Big) \oplus \{N\}'  
=\sub{\cE}{N}^s\oplus \sub{T}{N}\orbu{D,\,\eps} \ ,
\end{equation}
where the subspaces $\sub{\cE}{N}^s$ are the same as those constructed in 
Theorem \ref{the key}. 
Since $D \in \diagdis$ then, with the notations of Theorem \ref{the key},  
$\displaystyle \rho = \max_{D'\in  \bola{D,\, \eps} } k_{D'} <1$  for  $\eps$ small enought.  Also, for every $N \in  \orbu{D,\,\eps}$, 
\begin{enumerate}
	\item Both $\sub{\cE}{N}^s$  and $\sub{T}{N}\orbu{D,\,\eps}$, are $T_N\,\Delta$-invariant;
	\item $\left\|\left. T_N\,\Delta\right|_{\sub{\cE}{N}^s}\right\|\leq\rho<1$,  and 
	$T_N\,\Delta \Big|_{\sub{T}{N}\orbu{D,\,\eps}} $ is the identity map of $\sub{T}{N}\orbu{D,\,\eps}$. 
	\end{enumerate}
The distribution of the  subspaces $\sub{\cE}{N}^s$ is still smooth, since the  
(oblique) projection  $E_N$ onto  $\sub{\cE}{N}^s$ parallel to $\sub{T}{N}\orbu{D,\,\eps}$ 
moves smoothly on $\orbu{D,\,\eps}$.
A brief justification of these facts can be found in the following Remark:

\begin{rem}\label{CFH}\rm 
Let $d = \frac{1-\rho}{3}\, $. Consider the open discs $\cU = \{z \in \C : |z|<\rho+d\}$ and 
$\cV = \{z \in \C : |1-z|<d\}$, which have disjoint closures. 
By Eq. \eqref{summma}, and items 1 and 2 of the previous discusion, one can deduce that 
the spectrum of $T_N\,\Delta$ is contained in $\cU\cup \cV$ for every $N \in \orbu{D,\,\eps}$. 
Moreover, 
if $f : \cU\cup \cV \to \C$ is the holomorphic map $f = \aleph_\cU\,$ (the charcateristic 
map of $\cU$), then 
$E_N = f(T_N\,\Delta )$ for every $N \in \orbu{D,\,\eps}$. If $\cM(\cU\cup \cV ) = 
\{ T \in \cM_{r^2}(\C) : \spec{T} \inc \cU\cup \cV \}$, which is an open 
subset of $\cM_{r^2}(\C)$, then the map 
$$
\cM(\cU\cup \cV ) \ni T \mapsto f(T) \peso{is of class $C^\infty$ } 
$$
(see Theorem 5.16 of Kato's book \cite{Kato}). 
Therefore, the distribution $\orbu{D,\,\eps} \ni N \mapsto E_N = f(T_N\,\Delta )$ is 
of class $C^\infty$. A similar type of argument can be used to show that 
$\orbu{D,\,\eps}$  is a smooth submanifold of $\mat$, for $\eps$ small enough. 
\EOE
\end{rem}

\begin{pro} \label{contdiagdis}  
The map 
$\aluit{\infty}{\cdot}$ is  of class $C^\infty$ on $\diagdis$, and 
the sequence $\{\aluit{n}{\cdot}\}_{n\in\N}\,$,  
resticted to  $\diagdis$,  converges uniformly on 
compact sets to  $\Delta^\infty (\cdot)$.
\end{pro}
\proof 
Let $T \in \diagdis$, denote 
$N = \aluit{\infty}{T}$ and let $D\in \diagdis$, a diagonal matrix such that $N\in \orbu{D}$. 
We can apply Theorem \ref{teorema 5.5} to the pair $\orbu{D,\,\eps} \inc \orb{D,\,\eps}$, for $\eps$ small. From now on, the proof follows the same steps as the proofs of 
Corollary \ref{entorno y algo mas}
 and Proposition \ref{retraccion}. 
\QED

\section{Proof of Theorem \ref{the key}}\label{prueba del teorema principal}\label{la prueba}

\subsection{Matricial characterization of $T_N \Delta$}
\noi Throughout this section we fix an invertible diagonal matrix $D\in\mat$ whose diagonal 
entries are denoted by $(d_1,\ldots,d_n)$. For every $j\in\{1,\ldots,n\}$, let  
$d_j=e^{\,i\theta_j}|d_j|$ be the polar decomposition of $d_j$, where $\theta_j\in [0,2\pi]$. 
Recall from Eq. \eqref{tan1} 
that the tangent space $T_D \orb{D}$ consists on those matrices 
$X \in \mat$ such that $X_{ij} = 0$ if $d_i = d_j\,$. 

\begin{fed}\label{producto de Hadamard}
Given  $A,B\in\mat$, $A\circ B$ denotes their Hadamard product, 
that is, if $A=(A_{ij})$ and $B=(B_{ij})$, then $(A\circ B)_{ij}=A_{ij} B_{ij}$. 
With respect to this product, each matrix $A\in\mat$ induces an operator $\had{A}$ 
on $\mat$ defined by $\had{A}(B)=A\circ B$, $B \in \mat$. 
\end{fed}

\begin{rem}\label{remark sobre hadamard}
Note that, by Eq. \eqref{tan1},  the subspace $\sub{T}{D}\orb{D}$ reduces the 
operator $\had{A}\,$, for every $A\in\mat$. 
This is the reason why, from now on, we shall consider all these operators as acting on 
$\sub{T}{D}\orb{D}$. Restricted in this way, it holds that
\[
\|\had{A}\|=\sup\{\|A\circ B\|_2: \ B\in\sub{T}{D}\orb{D} \ \mbox{and}\  
\|B\|_2=1\}=\max_{d_i\neq d_j} 
|A_{ij}| \ .
\]
\EOE
\end{rem}

\noi 
Let $\preal$ and $\pim$ be the projections defined on $\sub{T}{D}\orb{D}$ by
$$
\preal(B)=\frac{B+B^*}{2} \peso{and} \pim(B)=\frac{B-B^*}{2}.
$$
That is, $\preal$ (resp. $\pim$) is the restriction to $\sub{T}{D}\orb{D}$ of the  orthogonal projection onto the subspace of hermitian (resp. anti-hermitian) matrices.  
Observe that, for every $K \in \matah$ (i.e., such that $K^* = -K$) and $B\in\mat$ it holds that 
\begin {equation}\label{commu}
K\circ \preal(B)=\pim(K\circ B)  \peso{and} K\circ \pim(B)=\preal(K\circ B)\ .
\end{equation}
Denote by $\qd\,$ 
the orthogonal projection from $\sub{T}{D}\orb{D}$ onto $(\sub{T}{D}\orbu{D})^\bot$. 

\begin{lem}\label{pruyucciun}
Let $J,K\in\mat$ be the matrices defined by
\begin{align*}
K_{ij}=
\begin{cases}
|d_j-d_i|\sgn(j-i)&\mbox{if $d_i\neq d_j$}\\
0                 &\mbox{if $d_i = d_j$}
\end{cases}
\peso{and} 
J_{ij}=\begin{cases}
(d_j-d_i)K_{ij}^{-1}&\mbox{if $d_i\neq d_j$}\\
1                   &\mbox{if $d_i = d_j$}
\end{cases},
\end{align*}
for $1\leq i,j\leq r$. Then
\begin{enumerate}
	\item For every $A\in\mat$, $AD-DA=J\circ K\circ A$.
	\item  It holds that $\qd \ = \had{J}\pim\had{J}^{-1}$.
	\item If $H \in \matsa$ (i.e., if $H^*=H$), then 
	$\qd \had{H}=\had{H}\qd\,$.
\end{enumerate}
\end{lem}
\proof $\ $
\begin{enumerate}
	\item It is enough to note that $(J\circ K)_{ij}=d_j-d_i$ and $(AD-DA)_{ij}=(d_j-d_i)A_{ij}$.
	\item Since $|J_{ij}|=1$ for every $1\leq i,j\leq r$, the operator $\had{J}$ is unitary 
	in $(\mat,\|\cdot\|_2)$. Hence, $\had{J}\pim\had{J}^{-1}$ is an orthogonal projection. 
	Recall that
  $$
    \sub{T}{D}\orbu{D}=\{AD-DA:\ A\in\matah\}.
  $$
By Eq. \eqref{commu}, $\pim \had{K}= \had{K} \preal\,$. 
Then, given $X=AD-DA\in \sub{T}{D}\orbu{D}$,
  $$
    \had{J}\pim\had{J}^{-1}(X)=\had{J}\pim\had{J}^{-1}(\had{J}\had{K}A)= 
    \had{J}\pim \had{K}(A)=\had{J} \had{K} \preal(A)=0.
  $$
  So, $\sub{T}{D}\orbu{D}\subseteq \ker(\had{J}\pim\had{J}^{-1})$. But,    
  $\dim\sub{T}{D}\orbu{D}=\dim \ker(\had{J}\pim\had{J}^{-1})$. Therefore, we have that 
  $\qd=\had{J}\pim\had{J}^{-1}$.
  \item It is clear that $\had{H}\had{J}=\had{J}\had{H}$. On the other hand, since $H$ is 
  hermitian, $\had{H}$ also commutes with the projection $\pim\,$.
\QED
\end{enumerate}

\begin{rem}\label{la matriz}\rm
Let $N \in \orbu{D}$ and let $\qn$ be the orthogonal projection from $\sub{T}{N}\orb{D}$  onto $\big(\sub{T}{N}\orbu{D}\big)^\bot$. 
Then 
$\dersin{\Delta}{N}$ has the following $2\times 2$ 
matrix decomposition
\begin{equation}\label{derivada en dos por dos}
\dersin{\Delta}{N}=
\begin{pmatrix}
\sub{A}{1N} & 0\\
\sub{A}{2N} & I   
\end{pmatrix}  \barr{l} \qn\\ I- \qn \earr \ ,
\end{equation}
because $\dersin{\Delta}{N}$ acts as the identity on $\sub{T}{N}\orbu{D}$. 
The next Proposition gives a characterization of the significative parts 
$\sub{A}{1N} = \qn \big(\dersin{\Delta}{N}\big) \qn$ and 
$\sub{A}{2N} = (I-\qn ) \big(\dersin{\Delta}{N}\big) \qn$ in the case $N=D$. \EOE
\end{rem}

\begin{pro}\label{uno}
Let  $\qd$ be the orthogonal projection 
onto $\big(\sub{T}{D}\orbu{D}\big)^\bot$. Then there exists $H \in \mat$ such that, 
if $H_1 = \preal (H)$ and $H_2 = \pim (H)$, 
$$
\qd \big(\dersin{\Delta}{D}\big) \qd =\qd\ \had{H_1}\ \qd \peso{and} (I-\qd) \big(\dersin{\Delta}{D}\big) \qd = (I-\qd)\ \had{H_2}\ \qd\, .
$$
The matrix $H_1$ can be characterized as   
\begin{align}
(H_1)_{ij}&=\frac{\big(1+e^{\,i(\theta_j-\theta_i)}\big)|d_i|^{1/2}|d_j|^{1/2}}{|d_i|+|d_j|} 
\peso {for every} 1 \le i, j \le r \ .
\label{h1} 
\end{align}
\end{pro}
\bdem

Fix a tangent vector $X=AD-DA \in T_D \orb{D}$, for some $A \in \mat$. Then
$$
\der{\Delta}{D}{X}=\left.\frac{d}{dt}\alu{e^{tA}De^{-tA}}\right|_{t=0}.
$$
Let $\gamma(t)=\big(e^{tA}De^{-tA}\big)^*\big(e^{tA}De^{-tA}\big)=e^{-tA^*}D^*e^{tA^*}e^{tA}De^{-tA}$. In terms of $\gamma$, we can write the curve $\alu{e^{tA}De^{-tA}}$ in the following way
$$
\alu{e^{tA}De^{-tA}}=\gamma^{1/4}(t)(e^{tA}De^{-tA})\gamma^{-1/4}(t).
$$
So, using that $(\gamma^{-1/4})'(0)=-\gamma^{-1/4}(0)\ (\gamma^{1/4})'(0)\ \gamma^{-1/4}(0)$ (which can be deduce from the identity $\gamma^{1/4}\gamma^{-1/4}=I$), we obtain
\begin{align*}
\der{\Delta}{D}{X}&= (\gamma^{1/4})'(0)\ D\gamma^{-1/4}(0)+\gamma^{1/4}(0)(AD-DA) \gamma^{-1/4}(0)\\&-\gamma^{1/4}(0)\ D\ \gamma^{-1/4}(0)\ (\gamma^{1/4})'(0)\ \gamma^{-1/4}(0)\\&=
(\gamma^{1/4})'(0)\ D|D|^{-1/2}+|D|^{1/2}(AD-DA) |D|^{-1/2}\\&-|D|^{1/2}\ D\ |D|^{-1/2}\ (\gamma^{1/4})'(0)\ |D|^{-1/2}\\
&=
\Big((\gamma^{1/4})'(0)\ D - D\ (\gamma^{1/4})'(0)\Big) |D|^{-1/2} +|D|^{1/2}(AD-DA) |D|^{-1/2}.
\end{align*}

\noi If we define the matrices $L,N\in\mat$ by
\begin{align*}
N_{ij}&=|d_j|^{-1/2},\\
L_{ij}&=|d_i|^{1/2}|d_j|^{-1/2},\\
\end{align*}
and take $J,K\in\mat$ as in Lemma \ref{pruyucciun}. Then
\begin{align*}
\der{\Delta}{D}{X}
&=N\circ( J\circ K \circ (\gamma^{1/4})'(0))+L\circ(J\circ K\circ A).
\end{align*}

\noi Now, we need to compute $(\gamma^{1/4})'(0)$. Firstly, we shall compute $(\gamma^{1/2})'(0)$, and then we shall repeat the procedure to get $(\gamma^{1/4})'(0)$. Using the identity $\gamma^{1/2}\gamma^{1/2}=\gamma$, we get
$$
\gamma^{1/2}(\gamma^{1/2})'+(\gamma^{1/2})'\gamma^{1/2}=\gamma'
$$
If $A=\gamma^{1/2}(0)$, $B=-\gamma^{1/2}(0)$ and $Y=\gamma'(0)$, we can rewrite the above 
identity in the following way
\[
A(\gamma^{1/2})'(0)-(\gamma^{1/2})'(0) B=Y.
\]
Therefore, $(\gamma^{1/2})'$ is the solution of Sylvester's equation $AX-XB=Y$. Using the well known formula for this solution (see \cite[Thm. VII.2.3]{[Bh]}), it holds that
\[
(\gamma^{1/2})'(0)=\int_0^\infty e^{-tA}Ye^{tB}\ dt= 
\int_0^\infty e^{-t\gamma^{1/2}(0)} \ \gamma '(0) \  e^{-t\gamma^{1/2}(0)}\ dt.
\]
In the same way, we get
\begin{align*}
(\gamma^{1/4})'(0)&=\int_0^\infty e^{-t\gamma^{1/4}(0)}\ (\gamma^{1/2})'(0)\ e^{-t\gamma^{1/4}(0)}\ dt\\ & 
= \int_0^\infty e^{-t\gamma^{1/4}(0)}\left(\int_0^\infty e^{-s\gamma^{1/2}(0)}\ \gamma'(0) \ 
e^{-s\gamma^{1/2}(0)}
\ ds\right)e^{-t\gamma^{1/4}(0)}\ dt \\ &
=\int_0^\infty\int_0^\infty e^{-\big(t\gamma^{1/4}(0) + s\gamma^{1/2}(0)\big) }\ \gamma'(0) \ 
e^{-\big(t\gamma^{1/4}(0)+s\gamma^{1/2}(0)\big)}\ 
ds\, dt.\\
\intertext{Finally, as $\gamma(0)=|D|^2$, we obtain}
(\gamma^{1/2})'(0) &
=\int_0^\infty\int_0^\infty e^{-(t|D|^{1/2} + s|D|)}\ \gamma'(0) \ 
e^{-(t|D|^{1/2}+s|D|)}\ 
ds\, dt.\\
\end{align*}
So, if $M\in\mat$ is the matrix defined by
\begin{align*}
M_{ij}&=\int_0^\infty\int_0^\infty e^{-(t|d_i|^{1/2}+s|d_i|)}\ 
e^{-(t|d_j|^{1/2}+s|d_j|)}\ ds\, dt\\
&=\int_0^\infty\int_0^\infty e^{-\Big(t(\ |d_i|^{1/2}+|d_j|^{1/2})\ +\ s(\ |d_i|+|d_j|)\Big)}\ 
\ ds\, dt\\
&=\int_0^\infty e^{-s\big(\ |d_i|+|d_j|\ \big)}\ \ ds\ 
\int_0^\infty e^{-t\big(\ |d_i|^{1/2}+|d_j|^{1/2}\big)} \ dt\\
&=\left.\frac{-e^{-s\big(\ |d_i|+|d_j|\ \big)}}{|d_i|+|d_j|}\right|_0^{\infty}\ \left.\frac{-e^{-t\big(\ |d_i|^{1/2}+|d_j|^{1/2}\big)}}{|d_i|^{1/2}+|d_j|^{1/2}}\right|_0^{\infty}\\
&=\frac{1}{|d_i|+|d_j|}\ \frac{1}{|d_i|^{1/2}+|d_j|^{1/2}},
\end{align*}
then $(\gamma^{1/4})'(0)=M\circ \gamma'(0)$. Our next step will be to compute $\gamma'(0)$. 
\begin{align*}
\gamma'(0)&=-A^*D^*D+D^*A^*D+D^*AD-D^*DA
           =2D^*\preal(A)D-(D^*DA+A^*D^*D)\\
          &=2D^*\preal(A)D-(D^*D\preal(A)+\preal(A)D^*D)-(D^*D\pim(A)-\pim(A)D^*D)
\end{align*}
Let $R,T^+,T^-\in\mat$ be the matrices defined by
$$
R_{ij}=2\bar{d}_i d_j \ , \quad 
T^+_{ij}=|d_i|^2+|d_j|^2 \ , \peso{and} 
T^-_{ij}=|d_j|^2-|d_i|^2 \ , \quad 1\le i,j \le r \ .
$$
Then, $\gamma'(0)$ can be rewritten in the following way
\begin{align*}
\gamma'(0)=R\circ\preal(A)-T^+\circ \preal(A)+T^-\circ \pim(A).
\end{align*}
In consequence, $\der{\Delta}{D}{AD-DA}$ can be characterized (in terms of $A$) as  
\begin{align*}
\der{\Delta}{D}{X} &=N\circ J\circ K\circ M\circ\Big[(R-T^+)\circ        
                       \preal(A)+T^-\circ\pim(A)\Big]+L\circ J\circ K\circ A.
\end{align*}

\noi Now, we shall express $\der{\Delta}{D}{X}$ in terms of $X=J\circ K\circ A$. 
Recall that, since  $K^*=-K$,  then  $\pim \had{K}= \had{K} \preal\,$, by Eq.  \eqref{commu}.
Therefore, 
\begin{align*}
\der{\Delta}{D}{X}&=M\circ N\circ (R-T^+)\circ J\circ \pim(K\circ A)\\
                      &+M\circ N\circ T^-\circ J\circ \preal(K\circ A)
                       +L\circ (J\circ K \circ A)\\
                      &=M\circ N\circ (R-T^+)\circ (\had{J}\pim\had{J}^{-1})(X)\\
                      &+M\circ N\circ T^-\circ (\had{J}\preal\had{J}^{-1})(X)
                       +L\circ (X) 
\end{align*}
Then, using Lemma \ref{pruyucciun}
\begin{align*}
\der{\Delta}{D}{X}&=\Big(M\circ N\circ (R-T^+)+L\Big)\circ \qd(X)\\
                      &+\Big(M\circ N\circ T^- +L\Big)\circ (I-\qd)(X).
\end{align*}
 
\noi If $H=M\circ N\circ (R-T^+)+L$, then $H_{i,j}=$ 
\begin{align*}
&=|d_i|^{1/2}|d_j|^{-1/2}+|d_j|^{-1/2}\frac{2\bar{d}_id_j-(|d_i|^2+|d_j|^2)} 
         {(|d_i|^{1/2}+|d_j|^{1/2})(|d_i|+|d_j|)}\\
       &=\frac{|d_i|^{1/2}|d_j|^{-1/2}(|d_i|^{1/2}+|d_j|^{1/2})(|d_i|+|d_j|)\ +\      2\bar{d}_id_j|d_j|^{-1/2}-|d_i|^2|d_j|^{-1/2}-|d_j|^{3/2}}{(|d_i|^{1/2}+|d_j|^{1/2})(|d_i|+|d_j|)}\\
       &=\frac{|d_i||d_j|^{1/2}+|d_i|^{3/2}+|d_i|^{1/2}|d_j|\ +\ 
       2\bar{d}_id_j|d_j|^{-1/2}-|d_j|^{3/2}}
       {(|d_i|^{1/2}+|d_j|^{1/2})(|d_i|+|d_j|)}\\
       &=\frac{|d_i||d_j|^{1/2}+|d_i|^{3/2}+|d_i|^{1/2}|d_j|+|d_j|^{3/2}\ +\ 
       2\bar{d}_id_j|d_j|^{-1/2}-2|d_j|^{3/2}}
       {(|d_i|^{1/2}+|d_j|^{1/2})(|d_i|+|d_j|)}\\
       &=1+2\frac{\bar{d}_id_j|d_j|^{-1/2}-|d_j|^{3/2}}  
       {(|d_i|^{1/2}+|d_j|^{1/2})(|d_i|+|d_j|)}.
\end{align*}
On the other hand
\begin{align*}
(M\circ N\circ T^- +L)&=|d_i|^{1/2}|d_j|^{-1/2}+|d_j|^{-1/2}\frac{|d_j|^{2}-|d_i|^{2}} 
                        {(|d_i|^{1/2}+|d_j|^{1/2})(|d_i|+|d_j|)}\\
                      &=|d_j|^{-1/2}\Big(|d_i|^{1/2}+|d_j|^{1/2}-|d_i|^{1/2}\Big)=1
\end{align*}
Therefore, we get that 
$
\der{\Delta}{D}{X}=\Big(\, H \qd +(I-\qd)\,\Big)(X).
$ 
Given $Y\in R(\qd)$,  
%
%
%
%
\begin{align*}
\qd\big(\dersin{\Delta}{D}\big) \qd (Y) &=\qd(H\circ Y)=(\had{J}\pim\had{J}^{-1})(H\circ Y)\\
&=J\circ\Big(\pim (H\circ \had{J}^{-1}Y)\Big)\\
&=\frac12 \ J\circ   \left( H\circ \had{J}^{-1}(Y)-\Big(H\circ \had{J}^{-1}(Y)\Big)^*\right) \\
&=\frac12 \ J\circ\Big( H\circ \had{J}^{-1}(Y)\,+\,H^*\circ \had{J}^{-1}(Y)\Big) \\
&=J\circ \preal(H)\circ \had{J}^{-1}(Y)= \preal(H)\circ Y = \qd \had{\preal(H)} ( Y ) \ .\\
\intertext{Analogously}
(I-\qd)\big(\dersin{\Delta}{D}\big) \qd (Y)  &=(I-\qd)(H\circ Y)=(\had{J}\preal\had{J}^{-1})(H\circ Y)\\
&=J\circ\Big(\preal (H\circ \had{J}^{-1}Y)\Big)\\
&=\frac12 \ J\circ\left(H\circ \had{J}^{-1}(Y)+\Big(H\circ \had{J}^{-1}(Y)\Big)^*\right)\\
&=\frac12 \ J\circ\Big( H\circ \had{J}^{-1}(Y)\,-\,H^*\circ \had{J}^{-1}(Y)\Big)\\
&=J\circ \pim(H)\circ \had{J}^{-1}(Y)= \pim(H) \circ Y = (I- \qd )  \had{\pim(H)} ( Y ) \ .
\end{align*}
So, Eq. \eqref{h1} holds. Moreover, 
\begin{align*}
(H_1)_{ij} &=\frac{1}{2}\left(1+2\frac{\bar{d}_id_j|d_j|^{-1/2}-|d_j|^{3/2}}                             {(|d_i|^{1/2}+|d_j|^{1/2})(|d_i|+|d_j|)} +     
                         1+2\frac{\bar{d}_id_j|d_i|^{-1/2}-|d_i|^{3/2}}       
                         {(|d_i|^{1/2}+|d_j|^{1/2})(|d_i|+|d_j|)}\right)\\
                        &=1+\frac{\bar{d}_id_j|d_j|^{-1/2}-|d_j|^{3/2} + 
                         \bar{d}_id_j|d_i|^{-1/2}-|d_i|^{3/2}}
                         {(|d_i|^{1/2}+|d_j|^{1/2})(|d_i|+|d_j|)}\\
                        &=\frac{|d_i||d_j|^{1/2}+|d_j||d_i|^{1/2} 
                         +\bar{d}_id_j|d_j|^{-1/2}+\bar{d}_id_j|d_i|^{-1/2}}
                         {(|d_i|^{1/2}+|d_j|^{1/2})(|d_i|+|d_j|)}\\
                        &=\frac{|d_i|^{1/2}|d_j|^{1/2}\Big(|d_i|^{1/2}+|d_j|^{1/2} 
                         +e^{\,i(\theta_j-\theta_i)}|d_i|^{1/2} 
                         +e^{\,i(\theta_j-\theta_i)}|d_j|^{1/2}\Big)}
                         {(|d_i|^{1/2}+|d_j|^{1/2})(|d_i|+|d_j|)}\\
                        &=\frac{\big(1+e^{\,i(\theta_j-\theta_i)}\big)
                         |d_i|^{1/2}|d_j|^{1/2}}{|d_i|+|d_j|}  \ , 
\end{align*}
which completes the proof.\edem


\begin{cor}\label{descomposicion}
Given $N \in \orbu{D}$, consider the matrix decomposition 
$$
\dersin{\Delta}{N}=
\begin{pmatrix}
\sub{A}{1N} & 0\\
\sub{A}{2N} & I   
\end{pmatrix}  \barr{l} \qn\\ I- \qn \earr \ ,
$$ 
as in Remark \ref{la matriz}. 
Then $\displaystyle \|\sub{A}{1N}\|\leq \max_{i,\,j\,:\ d_i\neq d_j}\frac{|1+e^{i(\theta_j-\theta_i)}|\, 
	|d_i|^{1/2}|d_j|^{1/2}}{|d_i|+|d_j|}<1$.
\end{cor}
\bdem
Let $N=UDU^*\in\orbu{D}$, for some $U \in \matu$. Then,
$$
\dersin{\Delta}{N}=\sub{Ad}{U}\Big(\dersin{\Delta}{D}\Big)\sub{Ad}{U}^{-1}
\peso{and} \sub{Q}{N}=\sub{Ad}{U}\Big(\sub{Q}{D}\Big)\sub{Ad}{U}^{-1}\ .
$$
Since $\sub{Ad}{U}:\sub{T}{D}\orb{D}\to\sub{T}{N}\orb{D}$ is an isometric isomorphism, it holds that
$$
\|\sub{A}{1N}\|=\left\|\qn \big( \dersin{\Delta}{N} \big) \qn \right\|=
\left\|\sub{Ad}{U}\Big(\qd\big( \dersin{\Delta}{D} \big)\qd \Big)\sub{Ad}{U}\inv \right\|=
\left\|\qd\big( \dersin{\Delta}{D} \big)\qd \right\|=\|\sub{A}{1D}\|.
$$ 
Take the selfadjoint matrix $H_1$ given by Proposition \ref{uno}.  Hence, 
$$
\displaystyle\|\sub{A}{1D}\|\leq \|\had{H_1}\|=\max_{i,\,j\,:\ d_i\neq d_j}\frac{|1+e^{i(\theta_j-\theta_i)}|\, |d_i|^{1/2}|d_j|^{1/2}}{|d_i|+|d_j|} \ .
$$
Finally, this maximum is strictly lower than one because, by the triangle inequality 
and the arithmetic-geometric inequality,  
	$$
	\frac{|1+e^{i(\theta_j-\theta_i)}|\, |d_i|^{1/2}|d_j|^{1/2}}{|d_i|+|d_j|}\ \leq \ 
	\frac{2\, |d_i|^{1/2}|d_j|^{1/2}}{|d_i|+|d_j|}\ \leq \ 1 \ .
	$$
But the equality holds  only if $\theta_j=\theta_i \mod(2\pi)$ and $|d_i|=|d_j|$, 
that is, if $d_i=d_j\,$.	
\edem
\begin{rem} It is easy to see, using Lemma \ref{pruyucciun} and Eq. \eqref{h1}, that 
$\dersin{\Delta}{D}$ is invertible, and therefore $\Delta$ is a local diffeomorphism
near $D$, if and only if $e^{i(\theta_j-\theta_i)}\neq -1$ for every $i,j$. 
This means that there are not pairs $d_i$, $d_j$ such that $d_i \cdot d_j \in \R_{<0}\,$.
\EOE
\end{rem}
\subsection{The proof}
Now we rewrite the statement of Theorem \ref{the key} and conclude its proof:
\begin{teo*}

The Aluthge transform $\alu{\cdot}:\orb{D}\to\orb{D}$ is a $C^\infty$ map, and for every $N\in\orbu{D}$, there exists a subspace $\sub{\cE}{N}^s$ in the tangent space $\sub{T}{N}\orb{D}$ such that
\begin{enumerate}
	\item $\sub{T}{N}\orb{D}=\sub{\cE}{N}^s\oplus \sub{T}{N}\orbu{D}$;
	\item Both, $\sub{\cE}{N}^s$ and  $\sub{T}{N}\orbu{D}$, are $T_N\,\Delta$-invariant;
	\item $\left\|\left. T_N\,\Delta\right|_{\sub{\cE}{N}^s}\right\|\leq\sub{k}{D}<1$, where $\displaystyle \sub{k}{D}=\max_{i,\,j\,:\ d_i\neq d_j}
	\frac{\big |1+e^{i(\arg(d_j)-\arg(d_i))}\big|\, |d_i|^{1/2}|d_j|^{1/2}}{|d_i|+|d_j|}$;
  \item If $U\in\matu$ satisfies $N=UDU^*$, then $\sub{\cE}{N}^s=U(\sub{\cE}{D}^s)U^*$.
  \end{enumerate}
  In particular, the map $\orbu{D} \ni N \mapsto \sub{\cE}{N}^s$ is smooth. This fact can be 
  formulated in terms of the projections $P_N $ onto $\sub{\cE}{N}^s$ parallel to $\sub{T}{N}\orbu{D}$, 
  $N \in \orbu{D}$.
  \end{teo*}

\bdem 
Fix $N=UDU^*\in\orbu{D}$. By Corollary \ref{descomposicion} $\|\sub{A}{1N}\|<1$, so the operator $I-\sub{A}{1N}$ acting on $R(\qn)$ is invertible. Let $\cE^s_N$ be the subspace defined by
$$
\cE^s_N=\left\{\begin{pmatrix}
y\\
-\sub{A}{2N}(I-\sub{A}{1N})^{-1} y
\end{pmatrix}:\ y\in R(\qn) \right\},
$$
where $\qn$, as in Corollary \ref{descomposicion}, is the orthogonal projection onto $\big(\sub{T}{N}\orbu{D}\big)^\bot$. A straightforward computation  shows that
$$
\sub{P}{N}=\begin{pmatrix}
I & 0\\
-\sub{A}{2N}(I-\sub{A}{1N})^{-1} & 0
\end{pmatrix}
\begin{array}{cc}
\qn\\
I-\qn
\end{array}\
$$
is a projection onto $\cE_N^s$ parallel to $\sub{T}{N}\orbu{D}$. Therefore
$$
\sub{T}{N}\orbu{D}=\cE^s_N \oplus \sub{T}{N}\orbu{D}.
$$
Moreover, since $\dersin{\Delta}{N}=\sub{Ad}{U}\Big(\dersin{\Delta}{D}\Big)\sub{Ad}{U}^{-1}$,
$\sub{Q}{N}=\sub{Ad}{U}\Big(\sub{Q}{D}\Big)\sub{Ad}{U}^{-1}$, and $\sub{P}{N}$ can be written as
$$\sub{P}{N}=\qn-(I-\qn)(\dersin{\Delta}{N})\qn\big(I-\qn(\dersin{\Delta}{N})\qn\big)^{-1}\qn,$$ it holds that 
$$\sub{P}{N}=\sub{Ad}{U}(\sub{P}{D})\sub{Ad}{U}^{-1}.$$
This shows that $\cE^s_N=U(\cE^s_D)U^*$ as we desired. On the other hand
\begin{align*}
\qn (\dersin{\Delta}{N})&=\begin{pmatrix}
\sub{A}{1N} & 0\\
\sub{A}{2N} & I
\end{pmatrix}
\begin{pmatrix}
I & 0\\
-\sub{A}{2N}(I-\sub{A}{1N})^{-1} & 0
\end{pmatrix}= \ 
\begin{pmatrix}
\sub{A}{1N} & 0\\
\sub{A}{2N}\Big(I-(I-\sub{A}{1N})^{-1}\Big) & 0
\end{pmatrix}\\&=
\begin{pmatrix}
\sub{A}{1N} & 0\\
\sub{A}{2N}\Big(-\sub{A}{1N}\Big)(I-\sub{A}{1N})^{-1} & 0
\end{pmatrix} \, \quad
= \ \begin{pmatrix} \sub{A}{1N} & 0\\
-\sub{A}{2N}(I-\sub{A}{1N})^{-1}\sub{A}{1N} & 0
\end{pmatrix}.\\
\intertext{and}
(\dersin{\Delta}{N})\qn & =\begin{pmatrix}
I & 0\\
-\sub{A}{2N}(I-\sub{A}{1N})^{-1} & 0
\end{pmatrix}
\begin{pmatrix}
\sub{A}{1N} & 0\\
\sub{A}{2N} & I
\end{pmatrix}= \ \begin{pmatrix}
\sub{A}{1N} & 0\\
-\sub{A}{2N}(I-\sub{A}{1N})^{-1}\sub{A}{1N} & 0
\end{pmatrix}.
\end{align*}
So, $\qn \dersin{\Delta}{N} = \dersin{\Delta}{N} \qn$. This implies that both, $\cE^s_N$ and $\sub{T}{N}\orbu{D}$, are invariant for $\dersin{\Delta}{N}$. Clearly, $\dersin{\Delta}{N}$ restricted to $\sub{T}{N}\orbu{D}$ is the identity. Hence, it only remains to prove that $\left.\big(\dersin{\Delta}{N}\big)\right|_{\cE^s_N}$ has norm lower or equal to $\sub{k}{D}$. Observe that it is enough to make the estimation at $\sub{T}{D}\orb{D}$. Indeed, for every $X\in\cE^s_N$, it holds that
$
\dersin{\Delta}{N}(X)=\sub{Ad}{U}\big(\dersin{\Delta}{D}\big)\sub{Ad}{U}^{-1}(X),
$
$\sub{Ad}{U}^{-1}(X)\in \cE^s_D$, and $\sub{Ad}{U}$ is an isometric isomorphism from $\sub{T}{D}\orb{D}$ onto $\sub{T}{N}\orb{D}$.  

\medskip
\noi So, let $Y=\begin{pmatrix}
y\\
-\sub{A}{2D}(I-\sub{A}{1D})^{-1} y
\end{pmatrix} \in \cE^s_D$. Then
\begin{align*}
\|(\dersin{\Delta}{D})\,(Y)\|_2^2&=
\left\|\begin{pmatrix}
\sub{A}{1D} & 0\\
\sub{A}{2D} & I
\end{pmatrix}\begin{pmatrix}
y\\
-\sub{A}{2D}(I-\sub{A}{1D})^{-1} y
\end{pmatrix}\right\|_2^2\\&=
\left\|\begin{pmatrix}
\sub{A}{1D}(y)\\
\sub{A}{2D}(y)-\sub{A}{2D}(I-\sub{A}{1D})^{-1} (y)
\end{pmatrix}\right\|_2^2\\&=
\left\|\sub{A}{1D}(y)\right\|_2^2 +\left\|\sub{A}{2D}(y)-\sub{A}{2D}(I-\sub{A}{1D})^{-1} (y)\right\|_2^2\\&\leq \sub{k}{D}^2
\left\|y\right\|_2^2 +\left\|-\sub{A}{2D}\sub{A}{1D}(I-\sub{A}{1D})^{-1} (y)\right\|_2^2.
\end{align*}
where the inequality holds because, by Corollary \ref{descomposicion}, $\|\sub{A}{1D}\|\leq\sub{k}{D}$.
On the other hand, by Lemma \ref{pruyucciun}, we know that $\had{H_1}\qd=\qd\had{H_1}\,$. So, 
 using  Proposition \ref{uno},  we obtain
\begin{align*}
\left\|-\sub{A}{2D}\sub{A}{1D}(I-\sub{A}{1D})^{-1} (y)\right\|_2^2&=
\left\|-(I-\qd)\,\had{H_2}\, \qd \, \had{H_1}\, \qd 
\Big((I-\sub{A}{1D})^{-1} (y)\Big)
\right\|_2^2\\&=
\left\|-\had{H_1}(I-\qd)\,\had{H_2}\, \qd \Big((I-\sub{A}{1D})^{-1} (y)\Big)\right\|_2^2\\&\leq 
\|\had{H_1}\|^2\left\|-(I-\qd)\,\had{H_2}\, \qd \Big((I-\sub{A}{1D})^{-1} (y)\Big)\right\|_2^2\\&=\sub{k}{D}^2\left\|-\sub{A}{2D}(I-\sub{A}{1D})^{-1} (y)\right\|_2^2.
\end{align*}
Therefore
\begin{align*}
\|(\dersin{\Delta}{D})\,(Y)\|_2^2&\leq \sub{k}{D}^2\ \|y\|_2^2+\sub{k}{D}^2\ \left\|-
\sub{A}{2D}(I-\sub{A}{1D})^{-1} (y)\right\|_2^2=\sub{k}{D}^2\|Y\|_2^2\ .
\end{align*}
The smoothness of the map $\orbu{D} \ni N \mapsto \sub{\cE}{N}^s$ 
follows from item (4) and the existence of $C^\infty$ local cross sections  
for the map $\pi_D : \matu \to \orbu{D}$, which exist by Proposition \ref{son variedades}. 
For example, if $\sigma_D : \cU \to \matu$ is such a section near $D$, then by 
item (4) and Eq. \eqref{ss} 
$$
P_N  = \sub{Ad}{\sigma_D (N)\,} P_D \sub{Ad}{\sigma_D (N)^* \,} \quad , \quad N \in \cU \ .
$$
This completes the proof.\edem

\appendix{
\section{Appendix: Stable manifold Theorem}

Let $f$ be a  smooth endomorphism of a Riemannian manifold and let $N$ be an f-invariant submanifold of $M$. Under the conditions of Theorem \ref{teorema 5.5}  we can suppose  that the tangent bundle at $N$ can be splitted in two $Df-$invariant subbundles, one given by the tangent bundle of $N$ and the other  being contracted by $Df$ (see Definition \ref{pseudohyperbolic}). In this case, as it holds for fixed points,  it is proved that for each point $x$ in $N$ there is a transversal smooth submanifold to $N$ containing 
$x$ and characterized by the points that converges  assymptoticaly to the orbit of $x$. The union of these submanifolds conforms a foliation in a neighborhood of $N$ (also called pre-lamination). This is the statement of theorem \ref{teorema 5.5}, which is obtained using a classical technique in dynamical systems  known as {\it graph transform operator} (see definition (\ref{gto})).
 This stable foliation has smooth leaves but in general is only continuous. However, if certain conditions over the $Df-$invariant  splitting   are also satisfied, then it can be proved that the foliation is smooth. This result, is consequences of the $C^r-$section theorem (stated here as theorem \ref{teorema 5.18} in subsection \ref{smooth}). Moreover, the $C^r-$section theorem  can be reformulated in a suitable version useful for our goals.  This  version is stated in theorem \ref{teorema 5.5bis}; in particular, in the statement is explicite which condition should be satisfied by the $Df-$invariant splitting (see inequality (\ref{r-acotacion})). To obtain this reformulation it is necessary to show that the graph transform operator introduced as a tool in the proof of the stable manifold theorem verifies certain properties. Therefore, and also for the sake of understanding for the reader, we give a sketch of the proof of the stable manifold theorem.

In our context, we want to apply the previous result for the case that the invariant submanifold is formed by fixed points. Therefore, we need to show that the hypothesis of theorem \ref{teorema 5.5bis} are full filed when we deal with a submanifold of fixed points. This is done in theorem \ref{teorema 5.5bisbis2}.

\subsection{Proof of theorem \ref{teorema 5.5}.}
\label{teorema 5.5ap}

\vskip 5pt

\noi
{\it Sketch of the proof:}
The proof consist in to use the graph transform operator. Basically consists in the following:  
In a neighborhood of any points $x\in N$ we consider the exponential map
$\exp_x: ({T_xM})_r\to M$ where $({T_xM})_r$ is the ball of radius $r$ in $T_xM$, and
we take the sets $$\hat{\cE^s}_x(r)= \exp(\cE^s_x\cap {(T_xM)}_r),\,\,\,\,\, \hat{\cF}_x(r)= \exp(\cF_x\cap {(T_xM)}_r).$$
Then it is taken $r$ small and the space of pre-lamination $\si$ such that for each $x\in N$ follows that $\si_x$ is a smooth map $\si_x: \hat{\cE^s}_x(r)\to \hat{\cF}_x(r)$ (in what follows, to avoid notation we simple note these subbundles with $\hat{\cE^s}_x$ and $\hat{\cF}_x$). 
Then it is taken the operator  which roughly speaking transform one pre-lamination into another one such that its images are related in the following way
(see (\ref{gto}) for details):
$$\si\to\tilde\si,\,\,\mbox{such that}\,\,\, image(\tilde \si_x) = f^{-1}(image(\si_{f(x)}))\cap B_r(x).$$
The goal is to prove that this operator is a contractive operator and so it has 
a fixed point. Latter it is shown that this fixed point corresponds to the stable 
lamination. 
%
%
Coming back to the sketch of the proof, first it is considered the maps 
$$ 
f^1_x= p^1_x\circ f: M\to \hat{\cE^s}_x \peso{and}  f^2_x= p^2_x\circ f: M\to \hat{\cF}_x \ ,
$$ 
where $p_x^1$ is the projection on $\hat{\cE^s}_x$ and $p_x^2$ is the projection on $\hat{\cF}_x.$
We take $$C^r(\hat{\cE^s}_x,\hat{\cF}_x)$$ the set of $C^r$ maps from $\hat{\cE^s}_x$ to $\hat{\cF}_x$ and we consider the space 
$$C^{r,0}(\hat{\cE^s},\hat{\cF})= \{\sigma: N\to C^r(\hat{\cE^s}_x,\hat{\cF}_x)\}$$ i.e.: for each $x\in N$ we take  $\si_x\in C^r(\hat{\cE^s}_x,\hat{\cF}_x)$ and we assume that $x\to \si_x$ moves continuously with $x$. We can represent $C^{r,0}(\hat{\cE^s},\hat{\cF})$ as a vector bundle over $N$ given by $N\times \{C^r(\hat{\cE^s}_x, \hat{\cF})\}_{x\in X}.$
Now we take the {\it graph transform operator}
\begin{eqnarray}\label{gto}
\Gamma_f(\sigma_x)= \big(f_x^2\circ(id,\si_{f(x)})\,\big)^{-1}
\circ \big(f_x^1\circ(id,\si_x)\,\big)\big | _{\hat{\cE^s}_x}.
\end{eqnarray}
It is proved that the graph transform operator is a contractive map and therefore it has a fixed point. In fact, to show that is contractive operator it is used the following remark

\begin{rem}\label{aclaraciones0}
	The Lipschitz constant of the graph transform operator is smaller that $\la$ where $\la$ 
	is the constant that bounds $\frac{||Df_{\cE^s}||}{m(Df_{|\cF})}$ (see inequality
	 (\ref{dominacion}) in Definition \ref{pseudohyperbolic}). In fact,  to prove that 
	it is enough to show that graph transform operator associated to $f$ is close to the 
	graph  transform operator $\Gamma_{Df}$  associated to $Df$and that $\la$ is an upper 
	bound for $Lip(\Gamma_{Df})$. The graph transform  operator associated to the derivative of $f$, acts on the space $L(\cE^s,\cF)$ which 
	is the bundle of linear maps from $\cE^s$ into $\cF$.  Using the splitting $\cE^s\oplus \cF$,  we can write $Df$ in the following way: 
	$$
	Df= \left [\begin{array}{clcr}
                          A &   0 \\
                          0  &  D   \end{array}\right ] \ ,
$$ 
where $A=Df_{|\cE^s}$ and $D=Df_{|\cF}\,$. Hence,  if $P\in L(\cE^s,\cF)\,$, then  $\Gamma_{Df}(P)$ is defined as \begin{eqnarray}
\label{lgto}
\Gamma_{Df}(P)=D^{-1}\circ P\circ A.
\end{eqnarray}
In particular, it follows that  $$Lip(\Gamma_{Df})= \frac{||A||}{m(D)}=\frac{||Df_{|\cE^s}||}{m(Df_{|\cF})}<\la<1.$$ Later, it is shown that the graph transform $\Gamma_f$ is close to $\Gamma_{Df}$ and so the remark follows. To see that `` $\Gamma_f$ is close to $\Gamma_{Df}$" observe that $D^{-1}$ in $x\in N$ is the derivative of ${f_x^2}^{-1}$ and  $A$ in $x\in N$ is the derivative of $f_x^1.$
\EOE
\end{rem}

\noi
 From the remark \ref{aclaraciones0}, we conclude that $\Gamma_f$ is a contractive operator with Lipschitz constant bounded by $\la.$

\vskip 5pt

\subsection{$C^r-$section theorem.} 
\label{smooth}

The goal is to prove that the pre-lamination obtained in Theorem \ref{teorema 5.5} is smooth. To do that, it is a used the following general theorem and latter we show how to adapt to prove the smoothness of the pre-lamination and we will address the particular case  of a submanifold of fixed points.

\begin{fed}
Let $\Pi: E\to X$ be a vector bundle with a metric space base $X$. We say that $d$ is an admissible on $E$ when:
\begin{enumerate}
\item it induces a norm on each fiber;

\item there is a Banach space $A$ such that the product metric on $X\times A$ induced $d$ on $E$;
\item the projection of $X\times A$ onto $E$ is of norm $1.$

\end{enumerate}
Without loss of generality we can assume that $E=X\times A$.

\end{fed}

\begin{fed}
Let $\Pi: E\to X$ be a vector bundle with a metric space base $X$, with an admissible metric on $E$. 
Let $X_0$ be a  subset  of $X$ and $D$  be the disc bundle of radius $C$ in $E$, where $C>0$ is a finite constant. Let $D_0$ be the restriction of $D$ to $X_0$; $D_0=D\cap \Pi^{-1}(X_0)$. Let $h$ be a continuous map  of $X_0$ into $X.$  We say that  $F:D_0\to D$ is a map which covers $h,$ if $$\Pi \circ F= h.$$

\end{fed}

\begin{teo}[$C^r-$section theorem.]\label{teorema 5.18}
Let $\Pi: E\to X$ be a vector bundle over the metric space $X$, with an admissible metric on $E.$ 
Let $X_0$ be a subset  of $X$ and $D$  be the disc bundle of radius $C$ in $E$, where $C>0$ is a finite constant. Let $D_0$ be the restriction of $D$ to $X_0$; $D_0=D\cap \Pi^{-1}(X_0)$. 
Let $h$ be an overflowing continuous map  of $X_0$ into $X$, that is $X_0\subset h(X_0).$ Let $F:D_0\to D$ be a map which covers $h.$ Suppose that there is  a constant $k,$ $0\leq k_0<1$ such that for all  $x\in X_0$, the restriction of $F$ to the fiber over $X$, $F_x:D_x\to D_{h(x)},$ is Lipschitz with constant at most $k$. Then:

\begin{enumerate}
\item There is a unique section $\sigma: X_0\to D_0$ such that $F(\mbox{Image of}\, \sigma)\cap D_0 = \mbox{Image of}\, \sigma.$
\item If, $X$, $X_0$ and $E$ are $C^r-$manifolds with bounded derivatives, if $\mu=Lip(h^{-1})$ be the Lipschitz constant of $h^{-1}$ and it is satisfied  \begin{eqnarray}\label{r-teo}
k\mu^r<1\end{eqnarray} then follows that $\sigma$ is $C^r$. \QED
\end{enumerate}
\end{teo}

\vskip 5pt

\noi
The previous theorem corresponds to theorem 5.18 of \cite{Sh1} (see page 58) and \cite{Sh2} (see page 44).

\vskip 5pt

\begin{rem}
Observe that in the previous Theorem, it is not assumed that the 
manifolds have to be compact. \EOE

\end{rem}

\subsection{Application to the smoothness of the stable lamination.}

\begin{teo}[Smoothness of the stable lamination]\label{teorema 5.5bis}
Let $f$ be a $C^r$ endomorphism of $M$ with a $\rho$-pseudo hyperbolic  submanifold $N$ with  $\rho< 1$. Let $\ewe^s: N\to \mbox{Emb}^r((-1,1)^k,M)$ be the  $C^r$-pre-lamination of class $C^0$,   introduced in Theorem \ref{teorema 5.5}.
If $m(\cdot)$ denotes the minimum norm, and 
\begin{eqnarray}\label{r-acotacion}
\frac{||Df_{/\cE^s}||}{m(Df_{/\cF})}||Df_{/\cF}||^r < \la <1,
\end{eqnarray} 
then $\ewe^s: \cU \cap N\to \mbox{Emb}^r((-1,1)^k,M)$ is a  $C^r$-pre-lamination of class $C^r$.

\end{teo}

\vskip 5pt
\noi
{\it Sketch of the proof:}
In the hypothesis of Theorem \ref{teorema 5.18} we consider $X=M$, $X_0=N$, $E=M\times \{C^r(\hat{\cE^s}_x, \hat{\cF}_x)\}_{x\in N}$ (i.e.: the pairs $(x,\si_x)$ such that $\si_x: \hat{\cE^s}_x\to \hat{\cF}_x),$ $h=f^{-1}$, 
$D_0= N\times \{C^r(\hat{\cE^s}_x, \hat{\cF})\}_{x\in N}$ and 
$F(x,\si)=(f(x), \Gamma_f) $ where $\Gamma_f$ is the graph transform operator associated to $f$. From remark \ref{aclaraciones0} follows that 
$Lip(F)$  is close to $\frac{||Df_{/\cE^s}||}{m(Df_{/\cF})}$  and it is immediate that $Lip(h^{-1})=Lip(f)= ||Df||.$ Therefore, if (\ref{r-acotacion}) holds, then 
$$
Lip(f)^rLip(\Gamma_f)<1 , 
$$ 
and therefore the inequality (\ref{r-teo}) holds and so we can apply  Theorem \ref{teorema 5.18}.
\QED

\vskip 5pt

\subsection{Application to a compact submanifold of fixed points.} 
Now we shows that we can apply \ref{teorema 5.5bis} to the case of a submanifold of fixed points.

\begin{cor}[Smoothness of the stable lamination for a submanifold of fixed points]\label{teorema 5.5bisbis2}
Let $f$, $M$ and $N$ as in Theorem \ref{teorema 5.5}. 
Let us assume that any point $p$ in $N$ is a {\bf fixed point}.  
Then $C^r$-pre-lamination $\ewe^s: \cN \to \mbox{Emb}^r((-1,1)^k,M)$ is  of class $C^r$.
\end{cor}

\begin{proof}
Observe  that $Df_{/\cF}=Id$. Therefore  
 $$
 \frac{||Df_{/\cE^s}||}{m(Df_{/\cF})}||Df_{/\cF}||^k=  ||Df_{/\cE^s}|| < \la <1
 $$
 and so it follows that  $\ewe^s: \cU \cap N\to \mbox{Emb}^r((-1,1)^k,M)$ is a  $C^r$-pre-lamination of class $C^r$, by Theorem \ref{teorema 5.5bis}.
\end{proof}

\begin{rem} \rm
Similar results to the one obtained in theorem \ref{teorema 5.5bisbis2} are obtained in \cite{PSW}. In this paper, it is shown that the stable foliation is $C^1$ assuming a similar condition  to (\ref{r-acotacion}) for the context of partial hyperbolic systems.
\end{rem}

}

\end{document}